\documentclass[a4paper,11pt]{article}

\usepackage[utf8]{inputenc}
\usepackage[english]{babel}
\usepackage{fancyhdr, todonotes}
\usepackage{indentfirst}

\usepackage{amssymb}
\usepackage{amsmath}
\usepackage{latexsym}
\usepackage{MnSymbol}
\usepackage{amsthm}
\usepackage{amscd}
\usepackage{amsmath}
\usepackage[all,cmtip]{xy}
\usepackage{multibib}

\usepackage[top=3cm, bottom=3cm, left=2.5cm, right=2.5cm]{geometry}

\newtheorem{thm}{Theorem}%[chapter]
\newtheorem{lem}[thm]{Lemma}%[chapter]
\newtheorem{prop}[thm]{Proposition}%[chapter]
\theoremstyle{definition}
\newtheorem{defn}{Definition}%[chapter]
\theoremstyle{remark}
\newtheorem{oss}[thm]{Remark}%[chapter]

\newcommand{\ZZ}{\mathbb Z}
\newcommand{\R}{\mathbb R}

\newcommand{\FF}{\mathbb F}
\newcommand{\NN}{\mathbb N}

\newcommand{\Qp}{\mathbb Q_p}
\newcommand{\OK}{\mathcal O_K}

\newcommand{\KK}{\overline K}

\newcommand{\DD}{\mathcal D}

\newcommand{\Ep}{\mathcal E_K(p)}
\newcommand{\Cp}{\mathcal C_K(p)}

\newcommand{\Gal}{\mathrm{Gal}}

\title{Upper ramification jumps in abelian extensions of exponent $p$}
\author{\textsc{Laura Capuano and Ilaria Del Corso}}
\date{}

\begin{document}

\maketitle

\begin{abstract}
In this paper we present a classification of the possible upper ramification jumps for an elementary abelian $p$-extension of a $p$-adic field.
The fundamental step for the proof of the main result is the computation 
of the ramification filtration for the maximal elementary abelian $p$-extension of the base field $K$. This result is a generalization of 
\cite[Lemma 9, p. 286]{Del_Corso_Dvornicich_2007}, where the same result is proved under the assumption 
that $K$ contains a primitive $p$-th root of unity. Using the class field theory and the explicit relations
between the normic group of an extension and its ramification jumps, it is fairly simple to recover necessary and sufficient conditions for
the upper ramification jumps of an elementary abelian $p$-extension of $K$.
\end{abstract}

\section{Introduction}

Let $K$ be a finite extension of $\Qp$. By the Hasse-Arf Theorem (\cite[p.76]{Serre}), the upper ramification jumps of a finite abelian 
extension $L/K$ are rational integers. The problem of determining the  upper ramification jumps sequence of an extension,
as well as the inverse problem to decide whether a set of integers could be the ramification jumps sequence of an extension with a fixed 
Galois group, is very difficult in general. However, this problem is completely solved in the case of cyclic extensions: 
a very neat result, due to Maus \cite{Maus_1973} in the case where $\zeta_p\nin K$, and to Miki 
\cite{Miki_1981} in the case where $\zeta_p\in K$, characterizes the sequence of integers which can be the ramification jumps of a totally 
ramified cyclic $p$-extension $L/K$. In this case, the ramification jumps completely determine 
the sequence of the ramification groups, since the quotients of the filtration are necessarily cyclic of order $p$. \\

In this paper we consider another basic case, namely the case of elementary abelian $p$-extensions of a $p$-adic field 
(some results for biquadratic extensions can be found in \cite{Byott_Elder_2002}). In this case, the 
ramification subgroups sequence depends upon the jumps and the order of the subgroups. 
In Theorem \ref{pea}, we characterize the sequences of couples of integers $(t,m)$, where $t$ denotes an upper 
jump and $m$ its ``size" (see Definition \ref{size}), which describes the ramification subgroups sequence of an elementary abelian 
$p$-extension of $K$. 

Our main tool is the class field theory and the explicit relation, already used in \cite{Sueyoshi_1984} for cyclic extensions,
between the normic group of an extension and its ramification jumps. The fundamental step for the proof of Theorem \ref{pea} is the computation 
of the ramification filtration for the maximal elementary abelian $p$-extension of the base field $K$. This result is contained in Theorem 
\ref{gr_ram} and is a generalization of \cite[Lemma 9, p. 286]{Del_Corso_Dvornicich_2007}, where the same result is proved under the assumption 
that $K$ contains a primitive $p$-th root of unity. \\

In the case of non-abelian extensions, the Hasse-Arf Theorem can fail and the upper ramification jumps can be not integral. However, one can give 
a classification for the lower ramification jumps. In this setting, very few cases are known. 
One special case can be found in \cite{Byott_Elder_2007} where, to better understand the counterexamples to the conclusion of the Hasse-Arf 
Theorem and as a first step towards an explicit description of wildly ramified Galois module structure, Byott and Elder classify 
the ramification break numbers of totally ramified quaternion extensions of dyadic number fields.

\section{Notation and preliminary results}

Throughout the paper $p$ will be a fixed prime number. If $K$ is a finite extension of $\Qp$, we shall denote by $e_K$ and $f_K$ the ramification 
index and the inertial degree of $K/\Qp$, and by $n_K$ the degree of the extension; hence, we have $n_K=e_Kf_K=[K:\Qp]$. 

We shall denote also by $\OK$ the ring of integers of $K$, by $\pi=\pi_K$ a uniformizer of $K$ 
(i.e. a generator of the maximal ideal $m_K$ of $\OK$) and by $v_K$ the valuation of $K$ normalized so that $v_K(\pi_K)=1$. We shall indicate 
the residue field of $K$ by $\KK$; then, $|\KK|=p^{f_K}$. 

Let $U_K$ be the group of unity of $K$, and consider its usual filtration $\{U_K^i\}_{i\ge1}$ given by  $U_K^i=1+m_K^i$ for $i\ge 1$.

For a finite extension $L/K$, we denote by $N_{L/K}$, $\DD_{L/K}$ and $\mathrm{Disc}_{L/K}$ the norm, the different and the discriminant of $L/K$ 
respectively. If $L/K$ is a Galois extension with Galois group $G$, we consider the filtration of $G$ given by the ramification subgroups:  
in our context, instead of the more classical lower numbering $G_i$ for the ramification subgroups, it is useful to use the upper numbering, 
so, for every $\nu \ge 0,$ we denote the ramification subgroups by $G^\nu$ (see \cite[Ch. IV]{Serre} for the definition and the fundamental properties of 
the ramification subgroups). \\

We recall here the following theorem which gives a rule to determine the ramification groups of a quotient (see \cite[Lemma 5, p. 75]{Serre}):
\begin{thm}[Herbrand] \label{Herbrand}
If $H$ is a normal subgroup of $G$, then, for every $\nu \gtr 0$,
$$ (G/H)^{\nu}=G^{\nu}H/H. $$
\end{thm}

We are interested in studying the filtration of the $G^{\nu}$ and, more specifically, the values of $\nu$ for which these subgroups change.
We give the following definition:

\begin{defn} We say that $s$ is a lower ramification jump for the extension $L/K$ if $G_s\neq G_u$ for every $s\gtr u$. Similarly, we say that $t$ is an 
upper ramification jump if $G^t\neq G^u$ for every $u\gtr t$.
\end{defn}

The lower jumps of an extension are always integers, whereas in general the upper jumps are not necessarily integers. However, in the case of abelian 
extensions, we have the following theorem (see \cite[p.76]{Serre}): 
\begin{thm}[Hasse-Arf]
If $G$ is an abelian group and if $\nu$ is a jump in the filtration $G^{\nu}$, then $\nu$ is an integer.
\end{thm}

As already observed in the introduction, the problem of determining whether a set of integers may be realized as the sequence of upper 
ramification jumps for an extension $L/K$ can be very difficult in general. 
The problem is completely solved in the case of cyclic extensions. \\

A necessary and sufficient condition for given $m$ natural numbers $t^1 < \cdots < t^m$ to be upper ramification jumps of a totally
ramified cyclic $p$-extension over $K$ was given by Maus (\cite{Maus_1973}, in two cases, namely when $\zeta_p \nin K$ and when, if $r$
is the maximal integer such that $\zeta_{p^r}\in K$, $v_K(\zeta_{p^r}-1)\nequiv 0 \mod p$), and by Miki \cite{Miki_1981} in the general case.
A constructive proof of the existence part of Miki's result was given by Sueyoshi in \cite{Sueyoshi_1984}. \\

The general result of Miki is rather technical to state; we at least recall what can happen in the easier case when $\zeta_p\nin K$:
% characterise the . Note that a cyclic extension of degree $p^m$ has exactly $m$ ramification jumps, as the 
%quotient between two consequent ramification groups has cardinality at most $p$. 
\begin{thm}[Maus, 1973] \label{Maus_thm}
Let $\{t^1 \less \cdots \less t^m\}$ be a finite set of integer numbers and suppose that $\zeta_p\nin K$; then, there exists a totally ramified cyclic
extension $L/K$ of degree $p^m$ with upper ramification jumps $t^1, \ldots, t^m$ if and only if the following conditions hold:
\begin{itemize}
 \item $1\le t^1 \less \frac{pe_K}{p-1}$ and $(t^1,p)=1$;
 \item if $t^i\less \frac{e_K}{p-1}$, then $t^{i+1}=pt^i$ or $pt^i\less t^{i+1} \less \frac{pe_K}{p-1}$ and $(p,t^{i+1})=1$;
 \item if $t^i\ge \frac{e_K}{p-1}$, then $t^{i+1}=t^i+e_K$.
\end{itemize}
\end{thm}
\ \\
Our aim is to characterize the upper ramification jumps and the ramification subgroups of an elementary abelian $p$-extension of $K$. 
In this case, ramification subgroups are clearly elementary abelian $p$-groups, so the ramification groups  sequence is completely determined by the jumps and the order of the subgroups. 

\begin{defn}\label{size}
Let $L/K$ be a Galois extension with Galois group $G$ and let $t\ge 1$ be an upper ramification jump. If $|G^t/G^{t+1}|=p^m$, we call $m$ the 
\textit{``size''} of the upper jump $t$.
\end{defn}

Given an elementary abelian $p$-extension, we can associate to the ramification subgroups a sequence of couples of integers $(t,m)$, where $t$ 
denotes an upper jump and $m$ its size. We will refer to the couple simply as to the ramification jump. \\

\bigskip
For convenience of the reader, we quote the class field correspondence theorem in a form which easily follows from Theorem 6.2 and 6.3 of 
\cite[Ch. III, p. 154]{Fesenko_Vostokov}.

\begin{thm} [Class field correspondence]\label{correspondence}
There is a one-to-one correspondence between the finite abelian extensions of $K$ and the subgroups of finite index of $K^{\times}$ given by
$L \longleftrightarrow N_{L/K}(L^{\times})$. This correspondence is an order reversing bijection between the lattice of finite abelian extensions of 
$K$ (with respect to the intersection $L_1 \cap L_2$ and the compositum $L_1L_2$) and the lattice of subgroups of finite index in $K^{\times}$
(with respect to the intersection $N_1 \cap N_2$ and the product $N_1N_2$). \\
Furthermore, if $L/K$ is the extension associated to the normic subgroup $N$ and $G$ is its Galois groups, $K^{\times}/N \cong G$, hence 
$|K^{\times}/N|=[L:K]$.
\end{thm}
\smallskip 
There is a strict connection between the ramification groups of an extension 
$L/K$ and the group $K^{\times}/N_{L/K}(L^{\times})$. In the case of totally ramified extensions of degree $p$, this is given by the following proposition:
%Let now assume that $L/K$ is a totally ramified extension. The following proposition holds:
%\begin{prop}
%The canonical map $U_K^n/N_{L/K}(U_L^{\psi(n)}) \longrightarrow K^{\times}/N_{L/K}(L^{\times})$ is injective, so the groups $U_K^n/N_{L/K}(U_L^{\psi(n)})$ form a decreasing filtration of $K^{\times}/N_{L/K}(L^{\times})$. We have $U_K^n/N_{L/K}(U_L^{\psi(n)})=0$ if and only if $G^n=\{1\}$. \\
%Furthermore let $c$ be the largest integer for which $G_c\neq \{1\}$ and let $f=\varphi(c)+1$. Then $U_K^f \subset N_{L/K}(L^{\times})$, and $f$ is the least integer enjoying
%this property.
%\end{prop}

%\noindent From this proposition we get the following useful theorem:
\begin{prop} \label{caracterization jump extension degree p}
Let $L/K$ be a totally ramified extension of degree $p$ and let $t$ be its upper ramification jump. Then
$$ t=\min\{j\in \NN\ |\ U_K^{j+1}\subseteq N_{L/K}(L^{\times}) \}. $$
\end{prop}

\begin{proof}
Let $t$ be the upper ramification jump of the extension $L/K$. 
Following \cite[Ch. IV]{Serre}, we denote by 
$\varphi_{L/K}: \R^+\to \R^+$ the Herbrand function %the function defined  by $$ \varphi_{L/K}(s)=\int_{0}^{s} {\cfrac{dx}{[G_0:G_x]}}, $$
and by $\psi_{L/K}$ its inverse. In this way, $G^{\nu}=: G_{\psi_{L/K}(\nu)}$ for every $\nu \gtr 0$.
In our case, an easy calculation gives
$$ \psi(n)=\left \{ \begin{matrix} n \quad \quad \quad \ \ \ \ \quad \mbox{ if $n\le t$} \\ \ \\ t+p(n-t) \quad \mbox{ \ if $n\gtr t$} \end{matrix} \right . $$
From \cite[Cor. 1, p. 228]{Serre}, we know that $N_{L/K}(U_L^{\psi(n)})\subseteq U_K^{n}$ for every $n$ and the equality holds if and only if $n\gtr t$. 
It follows that  $U_{K}^{t+1}\subset N_{L/K}(L^{\times})$, that is
$$ t\ge \min\{j\in \NN\ |\ U_K^{j+1}\subseteq N_{L/K}(L^{\times}) \}. $$
We have now to show that $t$ is exactly the minimum. If not, we would have $U_K^{t}\subseteq N_{L/K}(L^{\times})$. From \cite[Thm 1, p. 227]{Serre} 
we know that, for $n\ge 0$, the canonical map induced by inclusion and projection $U_K^n/N_{L_K}(U_L^{\psi(n)}) \longrightarrow K^{\times}/N_{L/K}(L^{\times})$ 
is injective, hence $U_K^n \cap N_{L/K}(L^{\times})=N_{L/K}(U_L^{\psi(n)})$.
This means that, if $U_K^n \subset N_{L/K}(L^{\times})$, we get $U_K^n \subset N_{L/K}(U_K^{\psi(n)})$, that is a contradiction if $t=n$. \\
\end{proof}

\section{The compositum of all extensions of degree $p$  of $K$} \label{compositum}

Let $ \Ep$ be the set of all the cyclic extensions of $K$ of degree $p$ within a fixed algebraic closure of $K$ and let $\Cp$ be the compositum 
of all extensions $E\in \Ep$; then, $\Cp$ is the maximal elementary abelian $p$-extension of $K$ in this fixed algebraic closure. In this section, 
following \cite{Del_Corso_Dvornicich_2007}, we determine the upper ramification jumps of $\Cp/K$.

\begin{prop} \label{max degree}
$$ [\Cp:K] = \left \{ \begin{matrix} n_K+1 \quad \quad \quad \ \ \ \mbox{ if $\zeta_p \nin K$} \\ \ \\ n_K+2 \quad \quad \quad \ \ \ \mbox{ if $\zeta_p \in K$} \end{matrix} \right . $$

\end{prop}
\begin{proof}
This is a classical result that can be easily proved using, for example, \cite[Ch. V, Prop 5.8 and Thm 5.7]{Narkiewicz}\footnote{This result
holds for every complete field with finite residue field.}.
\end{proof}

In \cite{Del_Corso_Dvornicich_2007}, the ramification subgroups of $\Cp/K$ are computed, using Kummer theory, in the case where $\zeta_p\in K$.  
The use of the class field theory allows us to generalize this result to a general field $K$. Also in this general case, the ramification 
groups can be computed via the study of all the subextensions of degree $p$. \\

Let $L/K$ be a Galois extension of degree $p$ and let $\DD_{L/K}=\pi_K^{D_L}$. Clearly $L\subset \Cp$ and, if $G_L=\Gal(\Cp/L)$, 
then $\Gal(L/K)\cong G/G_L$.

From the ramification-discriminant formula \cite[Prop. 4, p.64]{Serre},
$v_L(\DD_{L/K})= \sum_{i\ge 0} (|G_i|-1)$. 
In our case, $|G_i|=p$ if $0\le i\le t$ and $1$ otherwise, hence $v_L(\DD_{L/K})=(p-1)(t+1)$ and 
the jump of this extension is $t=\frac{D_L}{(p-1)}-1$. Hence,
$$(G/G_L)^{\nu}=(G/G_L)_{\nu}= \left \{ \begin{matrix} \ZZ/p \ZZ   & \text{if } \nu \le \cfrac{D_L}{(p-1)}-1 \\ \ \\
                                                  0 & \text{if } \nu\gtr \frac{D_L}{(p-1)}-1 \end{matrix} \right.    .$$
This information and Herbrand's Theorem \ref{Herbrand} allow us to reconstruct the ramification groups of $G$. In fact, by Herbrand's Theorem, 
$(G/G_L)^{\nu}=G^{\nu}G_L/G_L$, for every $\nu\ge 0$, hence
$$ (G/G_L)^{\nu}=0 \quad \iff G^{\nu}\subseteq G_L \iff \nu \gtr \frac{D_L}{(p-1)}-1. $$
Since $G_L$ runs over all subgroups of index $p$ of $G$ as $L$ runs over all normal extension of degree $p$ of $K$, it follows that
$$ G^{\nu}=\begin{LARGE}\bigcap \end{LARGE}_{\begin{small} \begin{matrix} L\subseteq \Cp \\ [L:K]=p \\ D_L/(p-1)\less \nu+1 \end{matrix} \end{small} } G_L. $$

\noindent This characterization of the ramification subgroups allow us to easily prove the following proposition: 

\begin{prop} \label{lemma subextension jump}
Let $M/K$ be a finite extension with Galois group $(\ZZ/p\ZZ)^h$; then $t$ is an upper ramification jump of $M/K$ if and only if there exists a 
subextension $L/K$ of degree $p$ with upper ramification jump equal to $t$.
\end{prop}

\begin{proof}
Assume that there exists a subextension $L\subseteq M$ such that $L/K$ is cyclic of degree $p$ with upper ramification jump $t$; we prove
that $t$ is an upper ramification jump for $M/K$. \\
Let $G$ be the group $\Gal(M/K)$ and $G_L=\Gal(M/L)$; then, $\Gal(L/K)\cong G/G_L=(G/G_L)^t$ and $(G/G_L)^{t+1}=\{1\}$ . \\
%  Suppose that $G^t=G^{t+1}$; by assumption $t$ is an upper ramification jump for $L/K$, hence $|(G/G_L)^{t+1}|=1$.
By Herbrand's Theorem \ref{Herbrand}, we have that, for every $s\ge 0$,
$$ (G/G_L)^{s}\cong G^{s}G_L/G_L \cong G^{s}/(G^{s}\cap G_L), $$
hence $(G/G_L)^{s}=\{1\}$ if and only if $G^{s}\subseteq G_L$. 
This proves that $G^{t+1}\subseteq G_L$, whereas $G^{t}\not\subseteq G_L$, hence $t$ is an upper ramification jump for $M/K$.\\

Assume now that $t$ is an upper ramification jump for $M/K$. From the previous description, we get
$$ G^t= \bigcap_{\begin{small} \begin{matrix} L\subseteq M \\ [L:K]=p \\ \frac{D_L}{(p-1)}\less t+1 \end{matrix} \end{small} } G_L
\quad \text{and} \quad G^{t+1}= \bigcap_{\begin{small} \begin{matrix} L\subseteq M \\ [L:K]=p \\ \frac{D_L}{(p-1)}\less t+2 \end{matrix} \end{small} } G_L. $$
Since $G^t\neq G^{t+1}$, there exists $L\subset M$ with $[L:M]=p$ and $D_L=(p-1)(t+1)$
and the upper ramification jump of this extension is exactly $t$.
\end{proof}

We want now to construct a normic group in $K^{\times}$ such that the corresponding extension is a subextension of $\Cp$ and 
its Galois group over $K$ has a given jump. 
To this aim, we need to describe the structure of the unit group $U_K.$   \\

Let  $I=\{ i \in \ZZ\ |\ 1\le i\less \frac{pe_K}{p-1} \mbox{ and } (p,i)=1 \}$, let $\KK$ be the residue field and let us fix a set 
$C=\{c_1, \ldots, c_{f_K} \}$ of elements
of $\OK$ such that the residues of its elements in $\KK$ form a basis of $\KK$ over $\FF_p$. If $\zeta_p\in K$, denote by $r$ the maximal
integer such that $K$ contains a $p^r$-root of unity.
%We have the following theorem on the structure of the group of units:
\begin{thm}[Fesenko-Vostokov, Ch. I, Prop. 6.4, p. 19] \label{FV}
Every $\alpha \in U_K^1$ can be written as a convergent product
$$ \alpha=\prod_{i\in I}\prod_{j=1}^{f_K} (1+c_j \pi^i)^{a_{ij}} {\omega_*}^{a}, $$
where:
\begin{itemize}
 \item if $\zeta_p \nin K$, $\omega_*=1$, $a=0$ and the above expression for $\alpha$ is unique, hence $U_K^1$ is a free $\ZZ_p$-module
of rank $n=e_Kf_K=[K:\Qp]$;
 \item if $\zeta_p \in K$, then $\omega_*=1+c_* \pi^{\frac{pe_K}{p-1}}$ is a principal unit such that $\omega_*\nin K^{p}$, $c_*\in C$
and $a\in \ZZ_p$. In this case, the above expression in not unique, and $U_K^1$ is a product of a free $\ZZ_p$-module of rank $n$ and the $p$-torsion
group $\mu_{p^r}$.
\end{itemize}
\end{thm}
Let us call
$ F=\{ (x,y)\in \ZZ \times \ZZ\ |\ x\in I,\ 1\le y\le f_K\} $; we put $\eta_{(x,y)}=1+c_y \pi^x$ for every $(x,y)\in F$.

It is known that the maximal unramified extension $K_{ur}$ of $\Cp$ is cyclic of degree $p$ (and it is the one associated to the
group $\langle {K^{\times}}^p, U_K^1 \rangle$).
% Using \cite[Ch. IV cor. p. 63]{Serre} we can reduce the study of the higher ramification
% groups (of index $i\ge 0$) to the totally ramified case.
The following lemma characterizes the maximal subextensions of $\Cp$ with only one ramification jump: 
\begin{lem} \label{lemma L_t} Let $t\in I$ and let $L_t/K$ the extension associated to the normic group
$$ N_t= \langle{K^{\times}}^p,\ \pi,\ \{ \eta_{(x,y)}\}_{(x,y)\in F,\ x\neq t},\ \omega_*  \rangle.
%\footnote{if $\zeta_p\nin K$ we put just $\omega_*=1$ } 
$$
Then, $L_t/K$ is an elementary abelian extension of degree $p^{f_K}$ with only one ramification jump equal to $t$.
\end{lem}

\begin{proof}
By Theorem \ref{correspondence}, we have $[L_t:K]=|K^{\times}/N_t|$. Now,  $K^{\times}/N_t \cong U_K^t/U_K^{t+1}$ and, since $t\in I$,
by Theorem \ref{FV} we have that $U_K^t/U_K^{t+1}\cong \overline K$. It follows that $[L_t:K]=|\overline K|=p^{f_K}$. 
Moreover, since $\langle {K^{\times}}^p, \pi \rangle \subset N_t$, then $L_t/K$ is totally ramified.

Using the previous proposition, it easy to see that $t$ is a ramification jump for $L_t$. In fact, we can consider the extension 
$L/K$ associated to the normic group
$$ N=\langle {K^{\times}}^p,\ \pi,\ \{ {\eta_{(x,y)}}_{(x,y)\in F \minus \{(t,1)\}}\}, \ \omega_* \rangle;$$
$L$ is a subextension of $L_t/K$ since $N_t \subset N$, has degree $p$ (because $|K^{\times}/N|=p$) and ramification jump equal to $t$ 
(this follows easily from Proposition \ref{caracterization jump extension degree p}). By Proposition \ref{lemma subextension jump}, 
we have that $t$ is also an upper ramification jump for $L_t/K$. 

We want to show that $t$ is the only possible jump. In fact, let $L$ be any sub\-extension of $L_t/K$ of degree $p$ over $K$; then, the group 
$N_L=N_{L/K}(L^{\times})$ is a subgroup of $K^{\times}$ of index $p$ and contains $N_t$, so $U_K^{t+1}\subset N_t \subset N_L$. 

On the other hand, 
$K^{\times}= \langle U_K^t, N_t \rangle $ and  $N_L\subsetneq K^{\times}$, so $U_K^t \nsubset N_L$ and, applying Proposition 
\ref{caracterization jump extension degree p}, 
its ramification jump is $t$. Using Proposition \ref{lemma subextension jump}, we get that $t$ is the only ramification jump of  $L_t/K$. \\
\end{proof}

\noindent If $\zeta_p\in K$, the field $\Cp$ has also a totally ramified subextension with jump not in the set $I$:
\begin{lem} \label{t'}
If $\zeta_p\in K$, let $t'=\frac{pe_K}{p-1}$ and let $L_{t'}$ be the extension associated to the normic subgroup
$$ N_{t'}= \langle {K^{\times}}^p, \pi,\ \{ \eta_{(x,y)}\}_{(x,y)\in F} \rangle. $$
Then, $L_{t'}/K$ is a cyclic extension of degree $p$ with ramification jump equal to $t'$.\footnote{We recall that, if $\zeta_p\in K$,
$(p-1) \ |\ e_K$ so ${\frac{pe_K}{p-1}}$ is an integer.}
\end{lem}
\begin{proof}
The argument of the proof is the same in the previous lemma. 
By the class field theory (Theorem \ref{correspondence}), we have that $[L_{t'}:K]=|K^{\times}/N_{t'}|$ and $|K^{\times}/N_{t'}|$ is a cyclic group generated 
by $\omega_* N_{t'}$ that has order $p$ (recall that $\omega_*$ is, by Theorem \ref{FV}, a principal unit of $U_K^{\frac{pe_K}{p-1}}$ such that 
$\omega_*\nin K^{p}$). The fact that the ramification jump is exactly $t'$ follows easily from Proposition 
\ref{caracterization jump extension degree p}, since $U_K^{t'}\nsubseteq N_{t'}$ 
(by construction $\omega_*\nin N_{t'}$) and $U_K^{t'+1}\subseteq {K^{\times}}^{p}\subseteq N_{t'}$. \\
\end{proof}
%\newpage
\begin{thm}\label{gr_ram}
If $\zeta_p \nin K$, the upper ramification groups of $\Cp/K$ are the following:
\begin{enumerate}
 \item $G=G^{-1}=(\ZZ/p\ZZ)^{n_K+1}$;
 \item $G^{0}=\ldots = G^{t^1}=(\ZZ/p\ZZ)^{n_K}$;
 %\item $G^{t^1+1}=\ldots=G^{t^2}=(\ZZ/p\ZZ)^{n_K-f}$;
 %\item[\vdots ] 
 \item $G^{t^{i}+1}=\ldots=G^{t^{i+1}}=(\ZZ/p\ZZ)^{n_K-if_K}$ for every $i=1, \ldots , e_K-1$;
 %\item[\vdots]
 %\item $G^{t^{e-1}+1}=\ldots= G^{t^e}=(\ZZ/p\ZZ)^{f}$;
 \item $G^{t^{e_K}+1}=\{e\}$;
\end{enumerate}
so, the upper ramification jumps are exactly $-1$ of size $1$ and $t^1 \ldots\ t^{e_K}$ of size $f_K$. \\
If $\zeta_p \in K$, the upper ramification groups of $\Cp/K$ are the following:
\begin{enumerate}
 \item $G=G^{-1}=(\ZZ/p\ZZ)^{n_K+2}$;
 \item $G^{0}=\ldots = G^{t^1}=(\ZZ/p\ZZ)^{n_K+1}$;
 %\item $G^{t^1+1}=\ldots=G^{t^2}=(\ZZ/p\ZZ)^{n_K+1-f}$;
 %\item[\vdots ] 
 \item $G^{t^{i}+1}=\ldots=G^{t^{i+1}}=(\ZZ/p\ZZ)^{n_K+1-if_K}$ for every $i=1, \ldots , e_K-1$;
 %\item[\vdots]
 %\item $G^{t^{e-1}+1}=\ldots= G^{t^e}=(\ZZ/p\ZZ)^{1+f}$;
 \item $G^{t^{e_K}+1}=G^{\frac{pe_K}{p-1}}=\{\ZZ/p\ZZ\}$;
 \item $G^{\frac{pe_K}{p-1}+1}=\{e\}$;
\end{enumerate}
so, the upper ramification jumps are exactly $-1$ and $\frac{pe_K}{p-1}$ of size $1$ and $t^1 \ldots\ t^{e_K}$ of size $f_K$.
\end{thm}

\begin{proof}
Let us consider the case $\zeta_p \in K$ (the case $\zeta_p \nin K$ is the same without taking into account the ``special subextension'' $L_{t'}$).
As done before, call $t'=\frac{pe_K}{p-1}$ and $I'=I \cup \{t'\}$.
Firstly, we show that $\Cp= K_{ur}L_{t'}\prod_{t\in I} L_t$. In fact, each extension on the right-hand side  is an elementary abelian 
$p$-extension, so $\Cp\supseteq K_{ur}\prod_{t\in I'} L_t$. On the other hand, $K_{ur}$ is linearly disjoint from  $\prod_{t\in I'} L_t$: 
in fact, for the class field theory and the previous constructions, the extension $K_{ur}/K$ is associa\-ted to the normic subgroup 
$\langle K^{{\times}p}, U_K^1 \rangle$, while $\prod_{t\in I'} L_t /K$ is associated to the normic subgroup $\bigcap_{t\in I'} N_t=\langle K^{{\times}p}, \pi \rangle $. 
Hence, the intersection of these extensions is the field associated to the normic group $ \langle K^{{\times}p}, \pi, U_K^1 \rangle =K^{\times}$, namely $K$. \\

With the same argument, one can show that, for every $\overline t\in I'$, the extension $L_{\overline t}$ is linearly disjoint from 
$\prod_{t\in I\setminus \{\overline t\}} L_t$. It follows that $$ [ K_{ur}\prod_{t\in I'} L_t:K]=[K_{ur}:K]\prod_{t\in I'} [L_t:K]=p^{n_K+2}=[\Cp :K],$$
so $\Cp= K_{ur}\prod_{t\in I'} L_t$. 

It is clear that all the integers $\{ -1,\ t^1, \ldots, \ t^{e_K},\ t' \}$ are upper ramification jumps for the extension $\Cp /K$, as all of them 
are upper ramification jumps for a cyclic subextension of $\Cp/K$ of degree $p$ (see the proofs of Lemma \ref{lemma L_t} and \ref{t'}). \\

To see that these are the only upper jumps, it is enough to prove that each of the jumps in $\{t^1,\ldots, t^{e_K}\}$ has at least size $p^{f_K}$ 
and $t'$ has size at least one. In this case, in fact, we get:
$$|G/G^0||G^{t'}/G^{t'+1}| \prod_{t\in I} |G^t/G^{t+1}| \ge p^{2+f_Ke_K}=|G|,$$
and this yields $|G^t/G^{t+1}|= p^{f_K}$ for every $t\in I$, $|G^{t'}/G^{t'+1}|=p$ and no more jump is possible.

We already know that $|G/G^0|=p$. Let $t\in I$ and call $H=\Gal(\Cp /L_t)$; by Galois correspondence, $G/H\cong \Gal(L_t/K)$.
From the previous lemma, we know that $L_t/K$ has only one upper ramification jump equal to $t$, so $(G/H)^t=(\ZZ/p\ZZ)^{f_K}$ and
$(G/H)^{t+1}=\{e\}$. 
On the other hand, Herbrand's Theorem \ref{Herbrand} ensures that, for each $s\gtr 0$, we have $(G/H)^{s} \cong G^{s}H/H$; moreover, 
$G^{s}H/H\cong G^{s}/G^{s}\cap H$, so we get $G^{t+1}\subseteq H$, and  
$$ |G^t/G^{t+1}|= |G^t/G^t\cap H|\cdot|G^t\cap H/G^{t+1}|=|(G/H)^t|\cdot|G^t\cap H/G^{t+1}|\ge p^{f_K}, $$ as wanted.
The same argument holds if we take $t'=\frac{pe_K}{p-1}$ and the extension $L_{t'}$ constructed in Lemma \ref{t'}.
\end{proof}

\section{The General Result} 

The following theorem classifies the sequence of couples of integers which can be the upper ramification jumps of an elementary
abelian $p$-extension of $K$.
\begin{thm}\label{pea}
Let $K$ be a finite extension of $\Qp$.
%Let $K$ be a complete discrete valuation field of $\mathrm{char}\ 0$ with a finite residue field of $\mathrm{char}\ p$.
Let $(t^1,m_1) \ldots (t^h,m_h)$ be couples of integers with $t^1\less \ldots \less t^h$; there exists an elementary abelian $p$-extension 
$M/K$ with upper ramification jumps $(t^1,m_1) \ldots (t^h,m_h)$ if and only if the following conditions hold:
\begin{enumerate}
   \item \label{1} for every $i=1, \ldots h,$ it holds $1\le t^i\less pe_K/(p-1)$ and $(t^i,p)=1$ with only two possible exceptions, 
         namely $t^1=-1$ and, in the case when $\zeta_p\in K$,  $t^h=\frac{pe_K}{p-1}$;
   \item \label{2} $1\le m_i\le f_K$ for every $i=1, \ldots h$, $m_1=1$  if  $t^1=-1$ and $m_h=1$ if $t^h=\frac{pe_K}{p-1}$.
\end{enumerate}
In this case, $[M:K]=\sum_{i=1}^hm_i$.
\end{thm}

\begin{proof}
Let $M$ be a subextension of $\Cp/K$ and  let $(t^1,m_1) \ldots (t^h,m_h)$ be its  ramification jumps. From Proposition 
\ref{lemma subextension jump}, we know that the jumps of $M/K$ are among those of $\Cp/K$, hence $t^1,\dots, t^h$ verify condition (\ref{1}).  
Moreover, denoting by $H$ the subgroup of $G={\Gal}(\Cp/K)$ fixing $M$, we have that the ramification filtration of $M/K$ is 
$(G/H)^{t^i}$ and $|(G/H)^{t^i}/(G/H)^{t^i+1}|=p^{m_i}$. Arguing as in the proof of Theorem \ref{gr_ram}, we have 
$$\left|\frac{(G/H)^{t^i}}{(G/H)^{t^i+1}}\right|= \left|\frac{G^{t^i}/G^{t^i+1}}{G^{t^i}\cap H/G^{t^i+1}\cap H}\right|=p^{f_K} /|G^{t^i}\cap H/G^{t^i+1}\cap H|,$$ 
hence $1\le m_i\le f_K$ if $t^i\in I$ and $m_1=1$ if $t^1=-1$. Moreover, 
if $\zeta_p\in K$ and $t^{h}=\frac{pe_K}{p-1}$, then $|G^{t^h}/G^{t^h+1}|=p$ as seen before, so $m_h=1$, 
namely, the $m_i$ verify condition (\ref{2}). \\ 
On the other hand, let $(t^1,m_1), \ldots ,\ (t^h,m_h)$ be a sequence verifying conditions (\ref{1}) and (\ref{2}); we construct an extension $M/K$ 
with these ramification jumps.  \\
For each $i= 1,\dots, h$ let $M_i/K$ be any subextension of degree $p^{m_i}$ of $L_{t^i}/K$, where $L_{t^i}$ is the extension defined in 
Lemma \ref{lemma L_t} when $t^i\in I$, $L_{t^h}=L_{t'}$ (the extension defined in Lemma \ref{t'}) if $t^h=\frac{pe_K}{p-1}$ and $L_{t^1}=K_{ur}$ if $t^1=-1$. 
Put $M=\prod_{i=1}^h M_i$; then, using the same techniques applied in Theorem \ref{gr_ram}, it is easy to see that the ramification jumps 
of $M/K$ are exactly $(t^1,m_1), \ldots,\ (t^h,m_h)$.
\end{proof}

\begin{oss}
If $\zeta_p\in K$, we can prove the same results using another approach, namely using Kummer theory (see \cite[p. 89]{Cassels_Frolich}). 
This method is more explicit as, 
in this way, we can construct an extension with fixed ramification jumps giving explicit generators. \\
In fact, by Kummer theory, we know that every $p$-extension is Galois and of the form $L=K(\sqrt[p]{a})$, with $a\in K^{\times}/ \setminus {K^{\times}}^p$. 
In particular, the element $a$ allows to determine  the ramification jump of the extension $L/K$. More precisely, 
the following proposition holds:
\begin{prop}
Let us take $a\in K^{{\times}}/{K^{{\times}}}^p$ and let us consider the extension $L=K(\sqrt[p]{a})$. Call $t$ the ramification jump; then:
\begin{itemize}
 \item if $0\less v_K(a)\less p$, then $t=\frac{pe_K}{p-1}$;
 \item if $v_K(a) = 0$ and $v_K(a-1)=l$ with $1\le l \less \frac{pe_K}{p-1}$ and $(l,p)=1$, then $t=\frac{pe_K}{p-1}-l$;
 \item if $v_K(a)=0$ and $v_K(a-1)=\frac{pe_K}{p-1}$, then $t=-1$ (and the extension is unramified).
\end{itemize}
\end{prop}

\begin{proof}
If $0\less v_K(a)\less p$ and $z$ is a $p$-th root of $a$, we can notice that, for every $h$ such that $(h,p)=1$, $K(z)=K(z^h)$, 
hence we can restrict to the case $v_K(a)=1$.

If $g$ is a generator of $\Gal(L/K)$, we have  
$v_K(g(\pi)-\pi)=v_L(z\pi-\pi)=v_L(z-1)+v_L(\pi)=\frac{pe_K}{p-1}+1,$ hence $t=\frac{pe_K}{p-1}$. 

If $a\in U_K^{l}\setminus U_K^{l+1}$, from the proof of \cite[Lemma 6, p. 15]{Del_Corso_Dvornicich_2007} we have that 
$v_L(\DD_{L/K})=(\frac{pe_K}{p-1}-l+1)(p-1)$. On the other hand, from the ramification-discriminant formula used in Section \ref{compositum}, 
$v_L(\DD_{L/K})=(p-1)(t+1)$. Comparing this with the previous expression, we have
$t=\frac{pe_K}{p-1}-l$. 
Finally, if $a\in U_K^{\frac{pe_K}{p-1}}$, using Hensel's Lemma \cite[p. 84]{Cassels_Frolich}, it is easy to see that the extension 
$L/K$ is unramified, hence the ramification jump is $t=-1$. 
\end{proof}

\noindent With this relation between the upper ramification jump and the valuation of the generator, one can easily give generators for the subextensions 
$L_{t}\subset \Cp$ constructed in Lemma \ref{lemma L_t}. Hence, the following proposition holds:
\begin{prop}
Take $t\in I$ and call $l=\frac{pe_K}{p-1}-t$; the extension $L_t=K\left (\sqrt[p]{ \eta_{(x,l)}},\ x=1, \ldots, f_K \right )$ 
is an elementary abelian extension of degree $p^{f_K}$ over $K$ with just one ramification jump equal to $t$. 

If $t=\frac{pe_K}{p-1}$, the extension $L_t=K(\sqrt[p]{\pi})$ is a totally ramified extension of degree $p$ with ramification jump equal to 
$\frac{pe_K}{p-1}$. 

If $t=-1$, the extension $L_{-1}=K(\sqrt[p]{\omega_*})$ is the only unramified extension of $K$ of degree $p$ and has ramification jump equal 
to $-1$.
\end{prop}
\end{oss}

\bigskip

\footnotesize
\begin{tabbing}
\hspace*{10cm} \= \hspace*{4cm} \kill
\parbox{\textwidth}{\textsc{Laura CAPUANO}} \>   \parbox{\textwidth}{\textsc{Ilaria DEL CORSO}} \\
Scuola Normale Superiore \> \parbox{\textwidth}{Dipartimento di Matematica} \\
Piazza dei Cavalieri, 7 \> \parbox{\textwidth}{Largo Bruno Pontecorvo, 5} \\
56127 Pisa (ITALY) \>  \parbox{\textwidth}{56127 Pisa (ITALY)} \\
E-mail: laura.capuano@sns.it \> \parbox{\textwidth}{E-mail:delcorso@dm.unipi.it}
\end{tabbing}

\bigskip

\footnotesize
\bibliographystyle{amsalpha} 
\bibliography{bibliography}

\end{document}